\documentclass{agtart_a}
\pdfoutput=1
\usepackage{pinlabel}

\title{The $C$--polynomial of a knot}

\author{Stavros Garoufalidis}
\givenname{Stavros}
\surname{Garoufalidis}
\address{School of Mathematics\\Georgia Institute of Technology\\\newline
Atlanta, GA 30332-0160\\USA}
\email{stavros@math.gatech.edu}
\urladdr{http://www.math.gatech.edu/~stavros}

\author{Xinyu Sun}
\givenname{Xinyu}
\surname{Sun}
\address{Department of Mathematics\\Mailstop 3368\\
Texas A\&M University\\\newline
College Station, TX 77843-3368\\USA}
\email{xsun@math.tamu.edu}
\urladdr{http://www.math.tamu.edu/~xsun/}

\volumenumber{6}
\issuenumber{}
\publicationyear{2006}
\papernumber{57}
\startpage{1623}
\endpage{1653}

\doi{}
\MR{}
\Zbl{}

\keyword{WZ algorithm}
\keyword{creative telescoping}
\keyword{colored Jones function}
\keyword{Gosper's algorithm}
\keyword{cyclotomic function}
\keyword{holonomic functions}
\keyword{characteristic varieties}
\keyword{$A$-polynomial}
\keyword{$C$-polynomial}
\subject{primary}{msc2000}{57N10}
\subject{secondary}{msc2000}{57M25}

\received{9 June 2005}
\revised{4 August 2006}
\accepted{29 August 2006}
\published{11 October 2006}
\publishedonline{11 October 2006}
\proposed{}
\seconded{}
\corresponding{}
\editor{CPR}
\version{}

\arxivreference{math.GT/0504305}

\AtBeginDocument{\let\bar\wbar\let\tilde\wtilde\let\hat\what}
\allowdisplaybreaks

\makeatletter
\def\cnewtheorem#1[#2]#3{\newtheorem{#1}{#3}[section]
\expandafter\let\csname c@#1\endcsname\c@proposition}
\makeatother

\theoremstyle{plain}
\newtheorem{theorem}{Theorem}

\cnewtheorem{lemma}[proposition]{Lemma}
\cnewtheorem{corollary}[proposition]{Corollary}
\theoremstyle{definition}
\cnewtheorem{definition}[proposition]{Definition}
\newtheorem{question}{Question}

\theoremstyle{remark}
\cnewtheorem{remark}[proposition]{Remark}

\def\BN{\mathbb N}
\def\BZ{\mathbb Z}

\def\BQ{\mathbb Q}

\def\BC{\mathbb C}

\def\calW{\mathcal W}

\def\A{\mathcal A}

\def\D{\Delta}

\def\a{\alpha}

\def\la{\langle}
\def\ra{\rangle}

\def\e{\epsilon}

\def\d{\delta}

\def\longto{\longrightarrow}

\def\Aloc{\A_{\mathrm{loc}}}
\def\Ch{\mathrm{Ch}}
\def\AJ{\mathrm{AJ}}
\def\SL{\mathrm{SL}}
\def\op{\mathrm{op}}

\def\Cert{\mathrm{Cert}}
\def\hatJ{\hat{J}}
\def\sgn{\mathrm{sgn}}

\def\rad{\mathrm{rad}}

\newcommand{\qfac}[1] {(q; q)_{#1}}
\newcommand{\qbinom}[2]{\binom{#1}{#2}_q}

\newcommand{\mat}[2][c|cccccccccccccccccccccccccccc]  
        {\left(\begin{array}{#1}#2\\\end{array}\right)}

\begin{document}

\begin{asciiabstract}
In an earlier paper the first author defined a non-commutative
A-polynomial for knots in 3-space, using the colored Jones
function. The idea is that the colored Jones function of a knot
satisfies a non-trivial linear q-difference equation. Said
differently, the colored Jones function of a knot is annihilated by a
non-zero ideal of the Weyl algebra which is generalted (after
localization) by the non-commutative A-polynomial of a knot.

In that paper, it was conjectured that this polynomial (which has to
do with representations of the quantum group U_q(SL_2)) specializes at
q=1 to the better known A-polynomial of a knot, which has to do with
genuine SL_2(C) representations of the knot complement.

Computing the non-commutative A-polynomial of a knot is a difficult
task which so far has been achieved for the two simplest knots. In the
present paper, we introduce the C-polynomial of a knot, along with its
non-commutative version, and give an explicit computation for all
twist knots. In a forthcoming paper, we will use this information to
compute the non-commutative A-polynomial of twist knots. Finally, we
formulate a number of conjectures relating the A, the C-polynomial and
the Alexander polynomial, all confirmed for the class of twist knots.
\end{asciiabstract}

\begin{htmlabstract}
<p class="noindent">
In an earlier paper the first author defined a non-commutative
A-polynomial for knots in 3-space, using the colored Jones function. The
idea is that the colored Jones function of a knot satisfies a non-trivial
linear q-difference equation. Said differently, the colored Jones
function of a knot is annihilated by a non-zero ideal of the Weyl algebra
which is generalted (after localization) by the non-commutative
A-polynomial of a knot.
</p>
<p class="noindent">
In that paper, it was conjectured that this polynomial (which has to do
with representations of the quantum group U<sub>q</sub>(sl<sub>2</sub>))
specializes at q=1 to the better known A-polynomial of a knot, which
has to do with genuine SL<sub>2</sub>(<b>C</b>) representations of the knot
complement.
</p>
<p class="noindent">
Computing the non-commutative A-polynomial of a knot is a difficult task
which so far has been achieved for the two simplest knots. In the present
paper, we introduce the C-polynomial of a knot, along with its
non-commutative version, and give an explicit computation for all
twist knots. In a forthcoming paper, we will use this information
to compute the non-commutative A-polynomial of twist knots. Finally,
we formulate a number of conjectures relating
the A, the C-polynomial and the Alexander polynomial, all confirmed
for the class of twist knots.
</p>
\end{htmlabstract}

\begin{abstract}
In an earlier paper the first author defined a non-commutative 
$A$--polynomial for knots in 3--space, using the colored Jones function. The 
idea is that the colored Jones function of a knot satisfies a non-trivial 
linear $q$--difference equation. Said differently, the colored Jones 
function of a knot is annihilated by a non-zero ideal of the Weyl algebra 
which is generalted (after localization) by the non-commutative 
$A$--polynomial of a knot.

In that paper, it was conjectured that this polynomial (which has to do
with representations of the quantum group $U_q(\mathfrak{sl}_2)$)
specializes at $q=1$ to the better known $A$--polynomial of a knot, which
has to do with genuine $\mathrm{SL}_2(\mathbb C)$ representations of the knot
complement.

Computing the non-commutative $A$--polynomial of a knot is a difficult task
which so far has been achieved for the two simplest knots. In the present
paper, we introduce the $C$--polynomial of a knot, along with its
non-commutative version, and give an explicit computation for all
twist knots. In a forthcoming paper, we will use this information
to compute the non-commutative $A$--polynomial of twist knots. Finally,
we formulate a number of conjectures relating
the $A$, the $C$--polynomial and the Alexander polynomial, all confirmed
for the class of twist knots.
\end{abstract}

\maketitle


\section{Introduction}
\label{sec.intro}

\subsection{The non-commutative $A$--polynomial of a knot}
\label{sub.nonA}

In \cite{Ga} the first author defined a non-commutative $A$--polynomial
for knots in 3--space, using the colored Jones function. The idea is that the
colored Jones function of a knot satisfies a non-trivial linear $q$--difference
equation. Said differently, the colored Jones function of a knot is annihilated
by a non-zero ideal of the Weyl algebra. By localizing, the Weyl algebra 
becomes a principal ideal domain, so that there is a single polynomial
generator, the non-commutative $A$--polynomial of a knot.

In \cite{Ga}, it was conjectured that this polynomial (which has to do
with representations of the quantum group $U_q(\mathfrak{sl}_2)$)
specializes at $q=1$ to the better known $A$--polynomial of a knot, which
has to do with genuine $\mathrm{SL}_2(\BC)$ representations of the knot
complement, Cooper--Culler--Gillet--Long--Shalen \cite{CCGLS}.

Computing the $A$--polynomial of a knot is a difficult task. For knots
with a small (about $10$) number of crossings, or with a small (about $7$)
number of ideal tetrahedra,
a numerical method was developed by Culler, see \cite{Cu}. 
For an alternative method that involves elimination, see Boyd \cite{Bo}.
For 2--bridge knots, simpler elimination methods are known.
All methods exhibit that the complexity of the $A$--polynomial 
(both with respect to the degrees of the 
monomials appearing, and with respect to their
coefficients) is exponential in the number of crossings.

\subsection{Can we compute the non-commutative $A$--polynomial?}
\label{sub.canwe}

At a first glance, it is not obvious that one can compute the
non-commutative $A$--polynomial of a knot. Let us explain a
theoretical algorithm for computation. Given a planar projection of a
knot with $c$ crossings, there is an explicit $c$--dimensional
multisum formula for the colored Jones polynomial, where the summand
is $q$--proper hypergeometric, see Garoufalidis and L{\^e}
\cite[Section 3]{GL1}. This has been implemented in Bar-Natan's {\tt
KnotAtlas} as a way of computing the colored Jones function of a knot,
see \cite{B-N}.

Given as input a multisum formula for the colored Jones polynomial,
the general theory of Zeilberger--Wilf computes a linear
$q$--difference equation by solving a system of linear equations; see
\cite{WZ}.  If one is lucky (and for general multisums unlucky cases
are known to exist) the linear $q$--difference equation is of minimal
order, thus computing the non-commutative $A$--polynomial. Even if one
is unlucky, there are costly factorization algorithms that in theory
will compute a minimal order $q$--difference equation; see
Petkov{\v{s}}ek, Wilf and Zeilberger \cite{PWZ}.

Using a computer implementation of the WZ method (Paule and Riese
\cite{PR3,PR2,PR1}), enabled the first author to give an explicit
formula for the non-commutative $A$--polynomial of the two simplest
knots: $3_1$ and $4_1$; see \cite{Ga}.

The main drawback of this implementation is that it works well when
the number of summation variables is $1$, but it becomes costly when the
number of summation variables increases. 

For 2--bridge knots, an alternative geometric method has been developed by 
Le that uses special properties of the Kauffman bracket skein module,    
L{\^e} \cite{Le}. Unfortunately, this method cannot be extended to the case
of non-2--bridge knots. In addition, the method is too costly to compute
the non-commutative $A$--polynomial of the $5_2$ knot.

Thus, two questions arise:

\begin{question}
\label{que.WZ1}
How can we reduce the number of summation variables in the WZ method?
\end{question}
\begin{question}
\label{que.521}
How can we compute the non-commutative $A$--polynomial of the $5_2$ and the
$6_1$ knots?
\end{question}

\subsection{The cyclotomic function of a knot}
\label{sub.jones}

To answer \fullref{que.WZ1}, 
we should look for efficient multisum formulas
for the colored Jones function of a knot. Thinking geometrically, it would be 
better to use a single variable for a whole sequence of twists between
two strands, rather than use one variable for each crossing.

As it turns out, Habiro \cite{H} introduced such formulas for the
colored Jones function of a knot.  It is a good moment to review the
colored Jones function, and Habiro's formulas.

A knot $K$ in 3--space 
is a smoothly embedded circle, considered up to 1--parameter
ambient motions of 3--space that avoid self-intersections. 
The {\em colored Jones function}
$J_K$ of a knot $K$ is a sequence of Laurent polynomials with integer
coefficients:
$$
J_K\co  \BN \longto \BZ[q^{\pm}].
$$
Technically, $J_K(n)$ is a quantum group invariant of the knot
colored by the $n$--dimen\-sional irreducible representation of $\mathfrak{sl}_2$,
normalized to $1$ for the unknot; see Turaev \cite{Tu}. When $n=2$, $J_K(2)$ is the
celebrated {\em Jones polynomial} of a knot, introduced in \cite{J}. 
One may think informally that the colored Jones function of a knot  encodes 
the Jones polynomial of a knot and its parallels.

In \cite{H}, Habiro introduced a key repackaging of the colored Jones
function $J_K$, namely the so-called {\em cyclotomic function} 
$$
\hatJ_K\co  \BN \longto \BZ[q^{\pm}],
$$
As the notation indicates, $\hatJ_K$ is in a sense a {\em linear
transformation} of $J_K$. More precisely, we have for every $n \geq 1$:
\begin{equation}
\label{eq.JC}
J_K(n)=\sum_{k=0}^\infty C(n,k) \hatJ_K(k)
\end{equation}
where
\begin{equation}
\label{eq.Cnk}
C(n,k)=\{n-k\} \dots \{n-1\} \{n+1\} \dots \{n+k\}
\end{equation}
and
$$
\{a\}=\frac{q^{a/2}-q^{-a/2}}{q^{1/2}-q^{-1/2}}.
$$
Notice that for every fixed $n$, the summation in Equation \eqref{eq.JC}
is finite, since $C(n,k)=0$ for $k \geq n$.

Habiro used {\em an integrality property} of the cyclotomic function
(namely, the fact that $\hatJ_K(n) \in \BZ[q^{\pm}]$ for all $n$) in order
to show that the Ohtsuki series of an integer homology sphere determines
its Witten--Reshetikhin--Turaev invariants, \cite{H}. The same integrality
property was used by Thang L\^e and the first author to settle the {\em
Volume Conjecture} to all orders, for small complex angles; see \cite{GL2}.

For our purposes, it is important that: 
\begin{itemize}
\item[(a)]
The transformation $J_K \longto \hatJ_K$ 
can be inverted to define $\hatJ_K$ in terms of $J_K$; 
see for example \cite[Section 4]{GL1}. 
\item[(b)]
$\hatJ_K$ satisfies a linear $q$--difference equation; see \cite{GL1}.
\item[(c)]
The cyclotomic function of twist knots has a single-sum formula; see
Equation \eqref{eq.Cpn} below.
\end{itemize}
 
We will use (c) to compute a minimal $q$--difference equation for the 
cyclotomic function of twist knots. Since $J_K$ determines and is
determined by $\hatJ_K$, in principle our results determine the 
non-commutative $A$--polynomial of twist knots. This motivates the results of 
our paper. En route, we will
introduce the $C$--polynomial of a knot and its non-commutative cousin.

Due to its length, the computation of the non-commutative
$A$--polynomial of twist knots will be postponed to a subsequent publication;
see \cite{GS}.

\subsection{What is a $q$--holonomic function and a $q$--difference equation?}
\label{sub.qholo}

Since we will be dealing with $q$--difference equations all along this paper,
let us review some general facts about the
combinatorics and geometry of $q$--difference equations.

There are two synonymous terms to $q$--difference equations: namely
{\em recursion relations}, and {\em operators}.
We will adopt the {\em operator point of view} when dealing with
recursion relations, in accordance to basic principles of 
physics and discrete math. An excellent reference is \cite{PWZ}.
Likewise, there is a synonymous term to a solution of a $q$--difference
linear equation: namely, a $q$--holonomic function.

For us, a {\em (discrete) function} $f$ is a map: 
$$
f\co \BN \longto \BQ(q)
$$
with values in the field of rational functions in $q$. Consider
two operators $E$ and $Q$ that act on the set of discrete functions by
$$
(Ef)(n)=f(n+1), \qquad (Qf)(n)=q^n f(n).
$$
It is easy to see that the operators $E$ and $Q$ satisfy 
$EQ=qQE$, and that $E$ and $Q$ generate a 
non-commutative {\em Weyl algebra} 
$$
\A=\BQ(q)\la Q, E \ra/(EQ-qQE).
$$
If $P=\sum_{j=0}^d a_j(Q,q) E^j$ is an element of $\A$, then the equation
$
Pf =0
$
is equivalent to the {\em linear $q$--difference equation}:
$$
\sum_{j=0}^d a_j(q^n,q) f(n+j)=0
$$
for all natural numbers $n$. Given a a discrete function $f$ as above, 
one may consider the set 
$$
I_f=\{P \in \A | Pf=0 \}
$$
of all linear $q$--difference equations that $f$ satisfies. 
It is easy to see that $I_f$
is a left ideal in $\A$. The following is a key definition:

\begin{definition}
\label{def.holonomic}
We say that $f$ is $q$--holonomic iff $I_f \neq 0$. 
\end{definition}
In other words, $f$ is $q$--holonomic iff it is a solution of a linear
$q$--difference equation. 
Unfortunately, the Weyl algebra $\A$ is not a principal (left)-ideal domain.
However, it becomes one after a suitable localization:
$$
\Aloc=\BQ(q,Q)\la E \ra/\left( E\a(Q,q)-\a(qQ,q)E \, | \, \a(Q,q) \in \BQ(Q,q)
\right).
$$
Moreover, the localized algebra still acts on discrete functions $f$.
Thus, given a $q$--holonomic function $f$, one may define
its {\em characteristic polynomial} $P_f \in \Aloc$, which is a generator
of the ideal $I_f$ over $\Aloc$. If we want to stress the dependence
of an operator $P$ on $E,Q$ and $q$, we will often write $P=P(E,Q,q)$.

There are three ways to view a $q$--holonomic function $f$:
\begin{itemize}
\item
The {\em $D$--module} $M_f=\Aloc/(P_f)$ and some of its elementary invariants:
its {\em rank}, and its {\em characteristic curve} 
$$
\Ch_f=\{(E,Q) \in \BC^2 | P_f(E,Q,1)=0 \}.
$$
The former is the $E$--degree of $P_f$ and the latter is a Lagrangian complex
curve in $\BC^2$.
\item
The {\em quantization} point of view, 
where we think of the operator $P_f(E,Q,q)$
as a $q$--deformation of the polynomial $P_f(E,Q,1)$. The zeros of the
latter polynomial define the characteristic curve, which is supposed
to be a classical object.  
\item
The {\em multi-graded} point of view. We may
think of $P_f(E,Q,q)$ as a polynomial in three
variables $E,Q$ and $q$ with integer coefficients. Then, $P_f(E,Q,q)$
and $P_f(E,Q,1)$ are, respectively, tri and bi-graded versions of $P_f(E,1,1)$.
\end{itemize}

Of course, $P_f(E,Q,q)$ is determined entirely by $f$.

As an example, consider the colored Jones function $J_K$ of a knot $K$,
and let $\AJ_K=\AJ_K(E,Q,q)$ denote its characteristic polynomial,
which here and below we will call the {\em non-commutative $A$--polynomial}
of the knot. 
The first author conjectured in \cite{Ga} that the evaluation of the
non-commutative $A$--polynomial at $q=1$ coincides
with the $A$--{\em polynomial} of a knot, whose zeros parametrize the  
 $\SL_2(\BC)$ character variety of the knot complement, restricted
to a boundary torus. For a definition of the $A$--polynomial, see
\cite{CCGLS}.

\subsection{The non-commutative $C$--polynomial of twist knots}
\label{sub.statement}

We now have all the ingredients to define the non-commutative
$C$--polynomial of a knot.

\begin{definition}
\label{def.C}
Given a knot $K$, let $C_K(E,Q,q)$ denote the characteristic polynomial
of its cyclotomic function $\hatJ_K$ and let $C_K(E,Q)$ denote $C_K(E,Q,1)$.
We will call $C_K(E,Q,q)$ (resp. $C_K(E,Q)$) the 
{\em non-commutative $C$--polynomial} (resp. the \em{$C$--polynomial}) of $K$.
\end{definition}

The reader should not confuse our $C$--polynomial with the 
{\em cusp polynomial} of a knot, due to X Zhang \cite{Zh}.

\begin{figure}[htpb]
\labellist\small
\pinlabel $p$ [l] at 35 110
\pinlabel full [l] <0pt, -10pt> at 35 110
\pinlabel twists [l] <0pt, -20pt> at 35 110
\endlabellist
\cl{\includegraphics[width=2in]{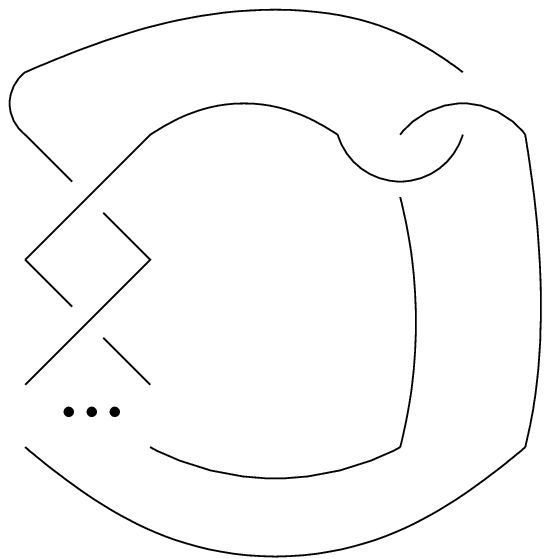}}
\caption{The twist knot $K_p$, for integers $p$}\label{twist}
\end{figure}

Consider the family of {\em twist knots} $K_p$ for integer $p$, shown 
in \fullref{twist}. The planar projection of $K_p$ 
has $2|p|+2$ crossings, $2|p|$ of which come from the full twists, 
and $2$ come from the negative {\em clasp}.

For small $p$, these knots may be identified with ones from Rolfsen's
table (see \cite{Rf}) as follows:
$$
K_1=3_1, \quad K_2=5_2, \quad K_3=7_2, \quad K_4=9_2
$$
$$
K_{-1}=4_1, \quad K_{-2}=6_1, \quad K_{-3}=8_1, \quad K_{-4}=10_1.
$$
Let $\hatJ_p(n)$ denote the cyclotomic function of $K_p$. Using
Masbaum \cite[Theorem 5.1]{Ma} (compare also with \cite[Section
3]{Ga}), it follows that:

\begin{equation}
\label{eq.Cpn}
\hatJ_p(n)= \sum_{k=0}^{\infty}     
        q^{n(n + 3)/2 + p k(k + 1) + k(k - 1)/2} (-1)^{n + k + 1}
        \frac{(q^{2k + 1} - 1)\qfac{n}} {\qfac{n + k + 1} \qfac{n-k}},
\end{equation}
where the {\em quantum factorial} and {\em quantum binomial coefficients}
are defined by:
$$
(A; q)_n = \left\{
\begin{tabular}{ll}
$(1-A) \cdots (1-Aq^{n-1})$ &   if $n>0$;  \\[1mm]
1       &       if $n = 0$;     \\
$\frac{1}{(1-Aq^{-1}) \cdots (1-Aq^{n})}$ & if $n < 0$, 
\end{tabular}
\right. 
$$
$$
\qbinom{m}{n} = \left\{
\begin{tabular}{ll}
$\frac{(q^{m-n+1};q)_n}{\qfac{n}}$ & if $ n \ge 0$; \\ [2mm]
$0$ & otherwise.                        
\end{tabular}
\right.
$$
We warn that we are using the unbalanced quantum factorials (common in discrete
math) and not the balanced ones (common in the representation theory of
quantum groups).

Equation \eqref{eq.Cpn} is the promised answer to \fullref{que.WZ1}
for the cyclotomic function of twist knots. For every fixed $p$, the summand
in \eqref{eq.Cpn} is $q$--{\em proper hypergeometric}
in the variables $n,k$. Notice that the summand  is {\em not} 
$q$--hypergeometric in all three variables $n,k,p$.

Our first result is an explicit formula for the non-commutative $C$--polynomial
of twist knots. 

\begin{definition}
\label{def.Cp}
\rm{(a)}\qua
For $p \in \BZ$, let us define $C_p(E,Q,q) \in \Aloc$ 
by:
\begin{eqnarray}
\label{eqn.recursion}
C_p(E,Q,q)=E^{|p|} + \sum_{i=0}^{|p|-1} a_p(Q, i) E^i,
\end{eqnarray}
where
\begin{equation}
\label{eq.apni}
a_p(q^n, i)=\left\{\hspace{-1mm}
\begin{tabular}{l}
$q^{(p-i) (n+p+1)} 
                \frac{\qfac{n+p-1}}{\qfac{n+i}} 
                \left( \sum_{j=0}^{i} q^{(2n+p + i + 1) j} 
                \qbinom{p-j}{p-i} \qbinom{p-i+j-1}{j}     \right.$ \\       
   \qquad$\left. - \sum_{j=0}^{i-1} q^{(2n+p+i+1) j + n+p} 
                \qbinom{p-j-1}{p-i} \qbinom{p-i+j-1}{j} \right)$          
        \qua  if $p>0$;       \\
0     \hspace{3.5in}  if $p = 0$;     \\
$q^{(p - i + 1)(n-p)}
                \frac{\qfac{n-p-1}}{\qfac{n+i}} 
                \left( -\sum_{j=0}^{i} q^{(2n-p+i) j}
                \qbinom{-p-j-1}{i-j} \qbinom{-p-i+j}{j}    \right.$\\
    \qquad$\left. + \sum_{j=0}^{i-1} q^{(2n-p+i)j+n-p}
                \qbinom{-p-j-2}{i-j-1} \qbinom{-p-i+j}{j} \right)    $       
\qquad if $p<0$.                       
\end{tabular}\hspace{-3mm} \right.
\end{equation}
In particular, $C_p(E,Q,q)$ is monic with respect to $E$ with
coefficients in $\BZ[Q^{\pm},q^{\pm}]$.

\rm{(b)}\qua
For $p \in \BZ$, let us define $C_p(E,Q) \in \BZ[E,Q^{\pm}]$ 
by:
\begin{equation}
\label{eq.CpEQ}
C_p(E,Q)=E^{|p|}+\sum_{i=0}^{|p|-1} b_p(Q,i) E^i,
\end{equation}
where
\begin{equation}
\label{eq.bpqi}
b_p(Q, i)  =  \left\{
\begin{tabular}{l}
$Q^{p-i} (1-Q)^{p-i-1}\left( \sum_{j=0}^{i} Q^{2j} 
\binom{p-j}{p-i} \binom{p-i+j-1}{j} \right. $ 
\\
\hspace{1in}$\left.- \sum_{j=0}^{i-1} Q^{2j+1} 
\binom{p-j-1}{p-i} \binom{p-i+j-1}{j} \right)$ \qquad\, if $p >0$; 
\\
0     \hspace{3.5in}  if $p = 0$;     \\
$Q^{-p-i+1} (1-Q)^{-p-i-1}\left( -\sum_{j=0}^{i} Q^{2j} 
\binom{-p-j-1}{i-j} \binom{-p-i+j}{j}\right. $
\\  
\hspace{1in}$\left.+ \sum_{j=0}^{i-1} Q^{2j+1} 
\binom{-p-j-2}{i-j-1} \binom{-p-i+j}{j} \right)$ \qquad if  $p <0$.
\end{tabular}
\right.
\end{equation}
\end{definition}

Notice that Equation \eqref{eq.bpqi} uniquely determines $a_p(Q,i)$.
A direct definition of $a_p(Q,i)$ would be cumbersome, since there is no
nice formula for $\binom{k}{l}_q$ as a rational function in $q$ and $Q=q^n$
when $k$ and $l$ are linear forms on $n$.

\begin{theorem}
\label{thm.1}
For every $p \in \BZ$, $C_p(E,Q,q)$ is the non-commutative $C$--polynomial
of the twist knot $K_p$.
\end{theorem}

An immediate corollary is:

\begin{corollary}
\label{cor.1}
For every $p \in \BZ$, $C_p(E,Q)$ is the $C$--polynomial of the twist knot 
$K_p$.
\end{corollary}
 
Our next result gives a 3--term recursion relation (with respect to $p$)
for the $C$--polynomial of twist knots.

\begin{theorem}
\label{thm.3}
{\rm(a)}\qua The $C$--polynomial of twist knots satisfies the 3--term relation:
\begin{equation}
\label{eq.3termp}
C_{p+2}(E,Q)=(Q-Q^2+E+Q^2E) \, C_{p+1}(E,Q)-Q^2 E^2 \,
C_{p}(E,Q)
\end{equation}
for all $p\geq 0$, with initial conditions:
$$
C_0(E,Q)=1, \qquad C_1(E,Q)=Q+E.
$$
Likewise, for $p \leq 0$ it satisfies a 3--term recursion
relation:
\begin{equation}
\label{eq.3termn}
C_{p-2}(E,Q)=(Q-Q^2+E+Q^2E) Q^{-2} \, C_{p-1}(E,Q)-E^2 Q^{-2}
\, C_{p}(E,Q),
\end{equation}
with intial conditions
$$
C_0(E,Q)=1, \qquad C_{-1}(E,Q)=-1+E.
$$
{\rm(b)}\qua Moreover,
\begin{equation}
\label{eq.CAlex}
C_p^{\op}(M-2+M^{-1},1)=\D_p(M)
\end{equation}
for all $p$.
\end{theorem}
Here, $\D_K(t) \in \BZ[t^{\pm}]$ denotes the {\em Alexander polynomial}
of a knot, normalized by $\D_K(t)=\D_K(t^{-1})$, and 
$\D_{\text{unknot}}(t)=1$; see \cite{Rf}.
Moreover, $\D_p(t)$ denotes the Alexander polynomial of the twist knot $K_p$.

In addition, if $P=P(E,Q)$ is a polynomial, then $P^{\op}$
\begin{equation}
\label{eq.Pop}
P^{\op}(E,Q)=E^{\deg_E P} P(E^{-1},Q)
\end{equation}
is essentially $P$, with its $E$--powers reversed.

\subsection{Relation between the $A$--polynomial and the $C$--polynomial
of twist knots}
\label{sub.relationAC}

The next theorem relates the $C$--polynomial 
of twist knots 
to the better-known $A$--polynomial of \cite{CCGLS}.
 In order to formulate our next theorem, we need to define a rational
map of degree $2$
\begin{equation}
\label{eq.phi}
\phi\co  \BQ(E,Q) \longto \BQ(L,M)
\end{equation}
by
$$
\phi(E)=\frac{L(M^2-1)^2}{M(L+M)(1+LM)}, \qquad
\phi(Q)=\frac{1+LM}{L+M}.
$$
For a motivation of this rational map, see \fullref{sub.motivation}.
\newpage

Let $A_p(L,M)$ denote the $A$--polynomial of the twist knot $K_p$.
The later has been computed by Hoste--Shanahan in \cite[Theorem 1]{HS1} (where 
is was denoted by $A_{J(2,2p)}$). It is known that the $A$--polynomial of a knot
in $S^3$ has even powers in $M^2$.
 
\begin{theorem}
\label{thm.CJp}
{\rm(a)}\qua For every $p \in \BZ$ we have: 
\begin{eqnarray}
\label{eq.Cphi}
\phi C_p^{\op}(E,Q) &=& A_p(L,M^{1/2}) \cdot 
\begin{cases}
\frac{(1+LM)^p}{M^p(L+M)^{3p-1}} & \text{if} \,\, p \geq 0; \\
\frac{1}{(M(L+M)(1+LM))^{|p|}} & \text{if} \,\, p <  0.
\end{cases}
\end{eqnarray}
{\rm(b)}\qua For every $p$, $C_p$ is an irreducible polynomial over $\BQ[E,Q]$.
\end{theorem}
\vspace{-5pt}

We can phrase the above theorem geometrically, as follows. The rational
map $\phi$ gives a rational map $\BC^2 \longto \BC^2$ where the domain
has coordinates $(E,Q)$ and the range has coordinates $(L,M)$. Then, we can
restrict the above map to the affine curves defined by $C_p^{\op}$ and $A_p$.
\vspace{-5pt}

\begin{corollary}
\label{cor.Cphi}
For all twist knots $K$, 
the map $\phi$ of \eqref{eq.phi} induces a Zariski dense map of degree $2$:
\begin{equation}
\label{eq.dominant}
\phi\co  \{ (L,M) \in \BC^2 \, | \, A_K(L,M)=0 \} \longto
\{ (E,Q) \in \BC^2 \, | \, C^{\op}_K(E,Q)=0 \}
\end{equation}
\end{corollary}
\vspace{-5pt}

Thus, one can associate two {\em plane curves} to a knot, namely the
$A$--curve, and the $C$--curve, which, in the case of twist knots, are
related by the map $\phi$ above. Thus one may consider their {\em
degrees} and their {\em genus}, discussed at length in Kirwan
\cite{Ki}.  The genus has the advantage of being a {\em birational
invariant}.

Rather than diverge to a lengthy algebraic geometry discussion, outside the
scope of the present paper, we state our next result here, and postpone
its proof in a subsequent publication.
 
\begin{theorem}
\label{thm.genus}
For every $p \in \BZ$, the genus of the $C_p(E,Q)$ polynomial is zero.
\end{theorem}
The proof uses the {\em Noether formula} for the genus of a plane curve
(see \cite[Theorem 7.37]{Ki}):
\begin{equation}
\label{eq.noether}
\text{genus}(C)=\frac{(d-1)(d-2)}{2} -\sum_P \delta(P)
\end{equation}
where $d$ is the degree and the (finite) sum is over the delta invariants
of the singular points of $C$. Since $\delta(P) > 0$ at the singular points
of $C$, and the left hand side is nonnegative, if one finds enough singular
points $P'$ such that the contribution makes the right hand side vanish,
then if follows that these are all the singular points of $C$ and moreover,
the genus of $C$ is zero. In our case, $d=3|p|-2$, 
the singular points $P'$ are 
$$
[0,0,1], \quad [1,0,0], \quad [0,1,1], \quad [0,1,0].
$$
in homogeneous coordinates, and their delta invariants are given by:
\begin{align*}
\d ([0,0,1]) &= |p|(|p|-1) &  \d([1,0,0]) &=(2|p|-3)(|p|-1) \\
\d([0,1,1]) &=|p|-1 & \d([0,1,0]) &=|p|-1.
\end{align*}
As a comparison, a Maple computation confirms that for 
$|p| \leq 30$ we have: 
\begin{equation}
\label{eq.genusAp}
\text{genus}(A_p(L,M^{1/2}))=
\begin{cases}
2p-2 & \text{if} \, p>0; \\
2|p|-1 & \text{if} \, p <0.
\end{cases}
\end{equation}
Unfortunately, the above method does not prove that the genus of the 
$A_p$ polynomial is given by \eqref{eq.genusAp} for all $p$, 
since it is hard to prove that the {\em only} singular points of $A_p$
for all $p$ are the ones suggested by Maple.

We thank N. Dunfield suggestions and for pointing the curious fact
about the genus of the $C_p$ polynomials.

\subsection{Plan of the proof}
\label{sub.plan}

As is obvious from a brief look, the paper tries to bring together two
largely disjoint areas: Quantum Topology and the Discrete Mathematics.
Thus, the proofs require some knowledge of both areas. We have tried to
separate the arguments in different sections, for different audiences.

In \fullref{sec.thmCJp}, we show that the sequence of polynomials
$C_p(E,Q)$ from Equation \eqref{eq.CpEQ} satisfy the 3--term recursion relations
\eqref{eq.3termp} and \eqref{eq.3termn}. Combining this result with
a 3--term recursion relation for the $A$--polynomial of twist knots (due
to Hoste--Shanahan), together with a matching of initial conditions allows
us to prove Equation \eqref{eq.Cphi}.
A side-bonus
of Equation \eqref{eq.Cphi} and of work of Hoste--Shanahan (that uses
ideas from hyperbolic geometry) is that the non-commutative polynomials
of Equation \eqref{eqn.recursion} are irreducible--
a property that the WZ algorithms cannot guarantee in general.

In \fullref{sec.WZ}, we give a crash course on the WZ algorithm
that computes recursion relations of sums of hypergeometric functions.
The ideas are beautiful and use elementary linear algebra. 
Using an explicit formula for the cyclotomic function of twist knots
(given in terms of a single sum of a $q$--hypergeometric function),
in \fullref{sec.twistproofs} we apply the WZ algorithm to confirm
that the cyclotomic function of twist knots satisfies the $q$--difference
equation 
$$
C_p(E,Q) \hatJ_p(n)=0
$$
for all $n \in \BN$. This, together with the  
irreducibility of $C_p(E,Q,q)$ obtained above, 
conclude the proof of \fullref{thm.1}.

In \fullref{sec.que} we present some
open questions (all confirmed for twist knots) about the structure of the
$C$--polynomial and the $A$--polynomial of knots.

Finally, in the Appendix we give a table of the non-commutative $C$--polynomial
of twist knots with at most $3$ crossings.

\subsection{Acknowledgements}
We wish to thank N Dunfield, J Rasmussen, D Zeilberger for many
stimulating conversations, and the anonymous referee for comments who
improved the presentation of the paper.

The first author was supported in part by National Science Foundation.

\section{Proof of Theorems \ref{thm.3} and \ref{thm.CJp}}
\label{sec.thmCJp}

\subsection[Proof of \ref{thm.3}]{Proof of \fullref{thm.3}}
\label{sub.thm3}

\begin{proof}(of \fullref{thm.3})
Consider the family of polynomials $C_p(E,Q)$ given by \eqref{eq.CpEQ}.
In this section we will show that this family of polynomials satisfies
the recursions stated in \fullref{thm.3}. Together with
\fullref{thm.1} (to be shown later), it will conclude the proof
of \fullref{thm.3}.

For convenience, we will convert the recursions in \eqref{eq.3termp}
and \eqref{eq.3termn} in 
backward shifts. That is, we define 
\begin{eqnarray*}
D_{p}(E, Q) &=& C_p^{\op}(E, Q) \\ &=& C_{p}(E^{-1},Q) E^{|p|}.
\end{eqnarray*}
Then we need to prove that
\begin{equation}
\label{eq.3term.backwardshift}
D_{p}(E,Q)=(QE - Q^2(E - 1) + 1) D_{p-1}(E,Q)-Q^2 D_{p-2}(E,Q).
\end{equation}
It  is clear from Equation \eqref{eqn.recursion.backwardshift} that
\begin{equation}
\label{eq.CpEQ.backwardshift}
D_p(E,Q)=1+\sum_{i=1}^{|p|} b'_p(Q,i) E^i
\end{equation}
and\eject
\begin{align*}
&b'_p(Q, i)  = Q^i (1-Q)^{i-1} 
\\
&\cdot\left\{ 
\begin{tabular}{l}
     $\left( \sum_{j=0}^{p-i} Q^{2j} 
\binom{p-j}{i} \binom{i+j-1}{j}     
 - \sum_{j=0}^{p-i-1} Q^{2j+1} 
\binom{p-j-1}{i} \binom{i+j-1}{j} \right)$\qquad\, if $p >0$; \\
$0$\hspace{4.1in} if $p=0$;\\
$Q^{2p+1}\left( -\sum_{j=0}^{-p-i} Q^{2j} 
\binom{-p-j-1}{i-1} \binom{i+j}{j}     
 + \sum_{j=0}^{-p-i-1} Q^{2j+1} 
\binom{-p-j-2}{i-1} \binom{i+j}{j}\right)$\\\hspace{4.2in}  if $p <0$.
\end{tabular}\right.
\end{align*}
When $p>0$, 
let 
$$s^{(1)}_p(E, Q, i, j) = Q^i(1 - Q)^{i - 1}Q^{2j}
\binom{p-j}{p-i-j}\binom{i+j-1}{j}E^i,$$
$$s^{(2)}_p(E, Q, i, j) = Q^i(1 - Q)^{i - 1}Q^{2j+1}
\binom{p-j-1}{p-i-j-1}\binom{i+j-1}{j}E^i,$$
$$D^{(1)}_p(E, Q) = \sum_{i=1}^{p}\sum^{p-i}_{j=0} s^{(1)}_p(E, Q),$$
$$D^{(2)}_p(E, Q) = \sum_{i=1}^{p}\sum^{p-i-1}_{j=0} s^{(2)}_p(E, Q),$$
and we have 
\begin{equation}
\label{eq.3term.summand}
D_p(E, Q) = 1 + D^{(1)}_p(E, Q) - D^{(2)}_p(E, Q).
\end{equation}
Using the same method mentioned before, 
it is easy to check that both $s^{(1)}_p(E, Q, i, j)$ and 
$s^{(2)}_p(E, Q, i, j)$ 
satisfy the {\em same} recursion
\begin{eqnarray*}
-Q^2 s^{(l)}_{p-2}(E, Q, i, j-1) - E (-1 + Q)Q s^{(l)}_{p-1}(E, Q, i-1, j) 
& & \\
 + Q^2 s^{(l)}_{p-1}(E, Q, i, j-1) + s^{(l)}_{p-1}(E, Q, i, j) - 
s^{(l)}_{p}(E, Q, i, j) &=& 0,
\end{eqnarray*}
for $l=1,2$.

Summing the above recursion over $i \ge 1$ and  $j \ge 0$, 
and noticing that $t^{(1)}_p(E, Q, i, j) = 0$ 
when $i > p$, $j < 0$, or $j > p-i$, 
we obtain
\begin{align}
\label{eq.DpEQ.1}
& -Q^2 D^{(1)}_{p-2}(E, Q) - E (-1 + Q)Q D^{(1)}_{p-1}(E, Q) + 
Q^2 D^{(1)}_{p-1}(E, Q)\\ 
&\hspace{2.5in}+ D^{(1)}_{p-1}(E, Q) - D^{(1)}_{p}(E, Q)\nonumber        \\ 
&= E (-1 + Q)Q s^{(1)}_{p-1}(E, Q, 0, 0) \nonumber      \\
&= - E Q.\nonumber 
\end{align}
Similarly
\begin{align}
\label{eq.DpEQ.2}
&  -Q^2 D^{(2)}_{p-2}(E, Q) - E (-1 + Q)Q D^{(2)}_{p-1}(E, Q) +
 Q^2 D^{(2)}_{p-1}(E, Q) \\
&\hspace{2.5in}+ D^{(2)}_{p-1}(E, Q) - D^{(2)}_{p}(E, Q)\nonumber   \\
&= E (-1 + Q)Q s^{(2)}_{p-1}(E, Q, 0, 0) \nonumber       \\
&= - E Q^2.\nonumber  
\end{align}
Now Equation \eqref{eq.3term.backwardshift} follows immediately
from Equations \eqref{eq.3term.summand}, \eqref{eq.DpEQ.1} and 
\eqref{eq.DpEQ.2},
which proves the theorem for $p > 0$.

For the case of $p < 0$, it is interesting that the 
backward-shifting 3--term recursion 
is the same as that when $p>0$.
To prove it, we only need to define, like before, 
$$s^{(1)}_p(E, Q, i, j) = Q^{2p+1+i} (1-Q)^{i-1} Q^{2j} 
\binom{p-j-1}{p-i-j} \binom{i+j}{j},$$
$$s^{(2)}_p(E, Q, i, j) = Q^{2p+1+i} (1-Q)^{i-1} Q^{2j+1} 
\binom{-p-j-2}{p-i-j-1} \binom{i+j}{j},$$
and realize that both of them satisfy the same recursion
\begin{eqnarray*}
-Q^2 s^{(l)}_{p-2}(E, Q, i, j) - E (-1 + Q) Q s^{(l)}_{p-1}(E, Q, i-1, j) 
& & \\
+ Q^2 s^{(l)}_{p-1}(E, Q, i, j-1) + s^{(l)}_{p-1}(E, Q, i, j) - 
s^{(l)}_{p}(E, Q, i, j-1) &=& 0,
\end{eqnarray*}
for $l=1,2$. This finishes the proof of part (a) of \fullref{thm.3}.

It remains to prove Equation \eqref{eq.CAlex}. The recursion relation for
$C_p$ from part (a) together with the initial conditions imply that
$$
C_p^{\op}(E,1)=1+pE.
$$
Since $\D_{K_p}(M)=1+p(M+M^{-1}-2)$, this concludes Equation 
\eqref{eq.CAlex} and \fullref{thm.3}.
\end{proof}

\subsection[Proof of \ref{thm.CJp}]{Proof of \fullref{thm.CJp}}
\label{sub.thm.CJp}

Consider the family of polynomials $C_p(E,Q)$ given by Equation 
\eqref{eq.CpEQ}. In \fullref{sub.thm3} we showed that $C_p(E,Q)$
satisfy the 3--term recursion relations \eqref{eq.3termp} and 
\eqref{eq.3termn}.

We will show that their evaluation $\phi C^{\op}(E,Q)$
satisfies Equation \eqref{eq.Cphi}.
In \cite[Theorem 1]{HS1}, Hoste--Shanahan give a 3--term recursion
relation for the $A$--polynomial of twist knots:
$$
A_p(L,M)=x(L,M) A_{p-\sgn(p)}(L,M) -y(L,M) A_{p-2\sgn(p)}(L,M),
$$
where
\begin{eqnarray*}
x(L,M) &=&
-L{+}L^2{+}2LM^2{+}M^4{+}2LM^4{+}L^2M^4{+}2LM^6{+}M^8{-}LM^8, \\
y(L,M) &=& M^4(L{+}M^2)^4.
\end{eqnarray*}
\fullref{thm.3} gives a 3--term relation for the $C$--polynomial of
twist knots. Assume, for simplicity, that $p >0$. Then, \fullref{thm.3}
implies that $C^{\op}_p(E,Q)$ satisfies a 3--term relation:
$$
C^{\op}_{p+2}(E,Q)=(EQ-EQ^2+1+Q^2) \, C^{\op}_{p+1}(E,Q)-Q^2 \,
C^{\op}_{p}(E,Q).
$$
Using the rational map $\phi$ of Equation \eqref{eq.phi}, a computation
shows that
$$
\phi(EQ-EQ^2+1+Q^2)=\frac{1+LM}{M(L+M)^2} x(L,M^{1/2}), \qquad
\phi(Q)=\frac{(1+LM)^2}{(L+M)^2}.
$$
Thus, it follows that both sides of \eqref{eq.Cphi} satisfy the same
3--term recursion relation for $p >0$.
Moreover, an explicit computation shows that \eqref{eq.Cphi} is verified
for $p=1,2$. The result follows for $p>0$, and similarly for $p <0$. 

For part (b), Hoste--Shanahan prove that the $A$--polynomial of twist
knots is irreducible; see \cite{HS2}. This, together with Equation 
\eqref{eq.Cphi} implies that any nontrivial factor of $C^{\op}_p(E,Q)$
must satisfy the property that its image under $\phi$ is a monomial
in $1+LM$ or $L+M$. This implies that any nontrivial factor of 
$C^{\op}_p(E,Q)$
will be of the form $(Q \pm 1)^2 + QE$. If $C^{\op}_p(E,Q)$ 
had any such factor, then evaluating at $Q=\mp 1$, 
it follows that $E$ divides $C_p(E,\mp 1)$. This is a contradiction,
by the explicit formula of \fullref{cor.1}. 
\qed

\section{A crash course on the WZ algorithm and Creative Telescoping}
\label{sec.WZ}

In this section we review briefly some key ideas of Zeilberger on recursion
relations of combinatorial sums. An excellent reference is \cite{PWZ},
which we urge the reader for references of the results in this section.
\vspace{5pt}

A {\em term} is $F(n, k)$ called {\em hypergeometric} 
if both $\frac{F(n+1, k)}{F(n, k)}$ and $\frac{F(n, k+1)}{F(n, k)}$
are rational functions over $n$ and $k$. In other words,
\begin{equation}
\label{eq.rationalF}
\frac{F(n+1, k)}{F(n, k)} \in \BQ(n,k), \qquad \frac{F(n, k+1)}{F(n, k)}
\in \BQ(n,k).
\end{equation}
Examples of hypergeometric terms
are $F(n,k)=(an+bk+c)!$ (for integers $a,b,c$), and ratios of products of such.
The latter are actually called {\em proper hypergeometric}.
A key problem is to construct recursion relations for sums of the form
\begin{equation}
\label{eq.sumS}
S(n)=\sum_k F(n,k),
\end{equation}
where $F(n,k)$ is a proper hypergeometric term.
The summation can be defined to be over all integers, even though in the 
cases that we consider in the paper, the summand vanishes for negative 
integers. Due to the telescoping nature of Sister Celine's method, we may
allow for a definite or indefinite summation range.
Sister Celine proved the following:
\vspace{5pt}

\begin{theorem}
Given a proper hypergeometric term $F(n, k)$, 
there exist an integer $I$ and a set of functions $a_i(n) \in \BQ(n)$, 
$0 \le i \le I$, such that
\begin{eqnarray} 
\label{eqn.celine}
\sum_{i=0}^{I} a_i(n) F(n+i, k) = 0.
\end{eqnarray}
\end{theorem}
\vspace{5pt}
The important part of the above theorem is that 
the functions $a_i(n)$ are independent of $k$. 
Therefore if we take the sum over $k$ on both sides, we get
\begin{eqnarray} 
\label{eqn.sum}
\sum_{i=0}^{I} a_i(n) \sum_{k}F(n+i, k) = 0.
\end{eqnarray}
In other words, we have:
\begin{equation}
\label{eq.haveS}
\sum_{i=0}^{I} a_i(n) S(n+i) = 0.
\end{equation}
So, Equation \eqref{eqn.celine} produces a recursion relation. How can we
find functions $a_i(n)$ that satisfy Equation~\eqref{eqn.celine}? The idea
is simple: divide Equation~\eqref{eqn.celine} by $F(n,k)$, and 
use \eqref{eq.rationalF} to convert the divided equation into an equation 
over the field $\BQ(n,k)$. Moreover,
$a_i(n)$ appear linearly. Clearing denominators, we arrive at an equation
(linear with respect to $a_i(n)$) over $\BQ(n)[k]$. Thus, the coefficients
of every power of $k$ must vanish, and this gives a linear system of
equations over $\BQ(n)$ with unknowns $a_i(n)$. If there are more unknowns
than equations, one is guaranteed to find a nonzero solution. By a counting
argument, one may see that if we choose $I$ high enough (this depends
on the complexity of the term $F(n,k)$), then we have more equations than
unknowns.
\vspace{5pt}

We should mention that although it can be numerically challenging to
find $a_i(n)$ that satisfy Equation~\eqref{eqn.celine}, it is routine
to check the equation once $a_i(n)$ are given. Indeed, one only need to divide 
the equation by $F(n,k)$,
and then check that a function in $\BQ(n,k)$ is identically zero.
The latter is computationally trivial.
\vspace{5pt}

This algorithm produces a recursion relation 
for $S(n)$. However, it is known that the algorithm
does not always yield a recursion relation of the smallest order.

Applying Gosper's algorithm, Wilf and Zeilberger invented another algorithm, 
the WZ algorithm. Instead of looking for $0$ on the right-hand side of 
Equation~\eqref{eqn.celine},
they instead looked for a function $G(n, k)$ such that
\begin{eqnarray} 
\label{eqn.wz}
\sum_{i=0}^{N} a_i(n) F(n+i, k) = G(n, k+1) - G(n, k).
\end{eqnarray}
Summing over $k$, and using telescoping cancellation of the terms
in the right hand side, we get a recursion relation for $S(n)$.
How to find the $a_i(n)$ and $G(n,k)$ that satisfy \eqref{eqn.wz}?
The idea is to look for a {\em rational function} $\Cert(n,k)$
(the so-called {\em certificate} of \eqref{eqn.wz}) such that
$$
G(n, k) = \Cert(n, k)F(n, k).
$$
Dividing out \eqref{eqn.wz} by $F(n,k)$ as before, one reduces this
to a problem of linear algebra. Just as before, given $a_i(n)$ and
$\Cert(n,k)$, it is routine to check whether \eqref{eqn.wz} holds.
\vspace{2pt}

Now, let us rephrase the above equations using operators.
Let us define two operators $N$ and $K$ that act on a function $F(n,k)$
by:
$$
(N F)(n, k) = F(n+1, k), \qquad 
(K F)(n, k) = F(n, k+1).
$$
Then, we can rewrite Equation~\eqref{eqn.wz} as
\begin{eqnarray} 
\label{eqn.operator}
\left(\sum_{i=0}^{I} a_i(n) N^i\right) F(n, k) = 
(K-1)G(n, k) = (K-1) \Cert(n, k) F(n, k).
\end{eqnarray}
Here, we think of $n$ and $k$ as operators acting on functions $F(n,k)$
by multiplication by $n$ and $k$ respectively. In other words,
$$
(nF)(n,k)=n F(n,k), \qquad (kF)(n,k)=k F(n,k).
$$
Beware that the operators $N$ and $n$ do not commute. Instead, we have:
$$
Nn=(n+1)N,
$$
and similarly for $k$ and $K$. 
\vspace{2pt}

Implementation of the algorithms are available in various platforms, 
such as,  Maple and Mathematica. See, for example, \cite{Z2} and \cite{PR2}.
\vspace{2pt}

Let us mention one more point regarding Creative Telescoping, namely
the issue of dealing with boundary terms. 
In the applications below, one considers
not quite the unrestricted sums of Equation \eqref{eq.sumS}, but 
rather restricted ones of the form:
\begin{equation}
\label{eq.sumSp}
S'(n)=\sum_{k=0}^\infty F(n,k),
\end{equation}
where $F(n,k)$ is a proper hypergeometric term.
When we apply the Creating Telescoping summation, we are left with some 
boundary terms $R(n) \in \BQ(n)$. In that case, Equation \eqref{eq.haveS}
becomes:
$$
\left(\sum_{i=0}^{I} a_i(n) N^i\right) S'(n) = R(n).
$$
This is an inhomogeneous equation of order $I$ which we can convert
into a homogeneous recursion of order $I+1$ by following trick:
apply the operator $$(N-1)\frac{1}{R(n)}$$ on both sides of the recursion, 
we get
$$
\left(\frac{1}{R(n+1)}N - 
\frac{1}{R(n)}\right)\left(\sum_{i=0}^{I} a_i(n) N^i\right) S'(n) = 0,
$$
$$\left(\frac{a_{I}(n+1)}{R(n+1)}N^{I+1}
  + \sum_{i=1}^{I}\left(\frac{a_{i-1}(n+1)}{R(n+1)} 
- \frac{a_i(n)}{R(n)} \right) N^i 
- \frac{a_0(n)}{R(n)}\right) S'(n) = 0.\leqno{\hbox{i.e.}}
$$
One final comment before we embark in the proof of the stated recursion
relations.
In Quantum Topology we are using $q$--factorials rather than factorials.
The previous results translate without conceptual difficulty to the $q$--world,
although the computer implementation costs more, in time.
A {\em term} is $F(n, k)$ called $q$--{\em hypergeometric} 
if 
$$
\frac{F(n+1, k)}{F(n, k)}, \frac{F(n, k+1)}{F(n, k)}
\in \BQ(q,q^n,q^k).
$$
Examples of $q$--hypergeometric terms are the quantum factorials of linear
forms in $n,k$, and ratios of products of quantum factorials
and $q$ raised to quadratic functions of $n$ and $k$.
The latter are called $q$--{\em proper hypergeometric}.

Sister Celine's algorithm and the WZ algorithm work equally well in the 
$q$--case. In either algorithms, we can (roughly speaking) replace $n$ and $k$ 
with $q^n$ and $q^k$ respectively, and the rest of the original proofs still 
apply naturally. 
The implementations of the $q$--case include \cite{PR3}, \cite{K} and \cite{Z2}.

\section{The non-commutative $C$--polynomial of twist knots}
\label{sec.twistproofs}

\subsection[Proof of \ref{thm.1}]{Proof of \fullref{thm.1}}
\label{sub.thm1}

First, let us make a remark for the trivial twist knot $K_0$.

\begin{remark}
\label{rem.p=0}
The colored Jones function of the trivial knot is $J_0(n)=1$ for all 
$n \geq 1$. Consequently, the cyclotomic function of the trivial knot is
$\hatJ_0(n)=\d_{n,0}$ (that is, $1$ when $n=0$ and $0$ otherwise).
The non-commutative $C$--polynomial of the trivial knot is $C_0(E,Q,q)=1$.
The $A$--polynomial of the trivial knot is $A_0(L,M)=1$ and the Alexander
polynomial of the trivial knot is $\D_0(M)=1$. This confirms all our theorems
for $p=0$.
\end{remark}

\begin{proof}(of \fullref{thm.1})
First we will prove that $\hatJ_p(n)$ satisfies the recursion 
relation:
\begin{equation}
\label{eq.newrec}
C_p(E,Q,q) \hatJ_p=0,
\end{equation}
where $C_p(E,Q,q)$ is given by Equation \eqref{eqn.recursion}.
We begin with rewriting the above equation as a recursion in 
backward shifts: 
\begin{eqnarray}
\label{eqn.recursion.backwardshift}
\hatJ_p(n) + \sum_{i=1}^{|p|} a'_p(n, i) \hatJ_p(n-i) = 0, 
\end{eqnarray}
where
$$
a'_p(n, i)=\left\{\begin{tabular}{l}
$       q^{i (n+1)} 
                \frac{\qfac{n-1}}{\qfac{n-i}}   
                \left( \sum_{j=0}^{p-i} q^{(2n - i + 1) j} 
                \qbinom{p-j}{i} \qbinom{i+j-1}{j}     \right.$        
\\
\hspace{0.56in}
                $\left. - \sum_{j=0}^{p-i-1} q^{(2n - i + 1) j + n} 
                \qbinom{p-j-1}{i} \qbinom{i+j-1}{j} \right)$\qquad
                if $p>0$;\\
0\hspace{3.6in} if $p = 0$;     \\
$q^{(2p + i + 1)n}
                \frac{\qfac{n-1}}{\qfac{n-i}}  
                \left( -\sum_{j=0}^{-p-i} q^{(2  n - i) j}
                \qbinom{-p-j-1}{i-1} \qbinom{i+j}{j}    \right.$  
\\      
\hspace{0.7in}
              $\left. + \sum_{j=0}^{-p-i-1} q^{(2n - i)j + n}                
              \qbinom{-p-j-2}{i-1} \qbinom{i+j}{j} \right)$\qquad       
if $p<0$.                       
\end{tabular} \right.
$$
When $p>0$, we define a number of functions for the purpose of convenience:
\begin{eqnarray*}
s_p(n, k)       &       =       &       
\frac{q^{n(n + 3)/2 + p k(k + 1) + k(k - 1)/2} 
(-1)^{n + k + 1}(q^{2k + 1} - 1)\qfac{n}}
{\qfac{n + k + 1} \qfac{n - k}},         \\[2.5mm]
t^{(1)}_p(n, k, i, j) &       =       &       
(-1)^i q^{-i(2n - i + 3)/2 + i(n + 1) + (2n - i + 1)j}  \\[2.5mm]
                        &               &
\frac {\qfac{n - 1} \qfac{n + k + 1} \qfac{n - k}}
{\qfac{n} \qfac{n + k - i + 1} \qfac{n - i - k}} 
\qbinom{p - j}{p-i-j}\qbinom{i + j - 1}{j},          \\[2.5mm]
t^{(2)}_p(n, k, i, j) &       =       &       
(-1)^i q^{-i(2n - i + 3)/2 + i(n + 1) + (2n - i + 1)j}  \\[2.5mm]
                        &  \cdot             &
\frac {\qfac{n - 1} \qfac{n + k + 1} \qfac{n - k}}
{\qfac{n} \qfac{n + k - i + 1} \qfac{n - i - k}} 
\qbinom{p - j - 1}{p-i-j-1}\qbinom{i + j - 1}{j},          \\[2.5mm]
r_p(n, k)       &       =       &       \sum_{i=1}^{p} 
\frac{a'_p(n, i) s_p(n-i, k)}{s_p(n, k)},           \\[2.5mm]
\Cert_p(n, k)       &       =       &       
\frac{q^{p k + p n + p} (q^{k + 1} - 1) (q^n - q^k)}
{(q^{2 k + 1} - 1) (q^n - 1)},                    \\[2.5mm]
D_p(n, k)       &       =       &       
\Cert_p(n, k) - \Cert_p(n, k-1) \frac{s_p(n, k-1)}{s_p(n, k)} - 1. 
\end{eqnarray*}
It is clear that 
$$
\sum_{k \ge 0}s_p(n, k) = \hatJ_p(n).
$$ 
Since 
$$
t^{(h)}_p(n, k, i, j) = 0  \quad \mathrm{if}\ j > p - i - h + 1
\ \mathrm{or}\ i>p,  \quad \mathrm{when}\ h = 1,2,
$$
and 
$$
\frac{s_p(n-i, k)}{s_p(n, k)} = 
(-1)^{-i} q^{-i(2n-i+3)/2} 
\frac{\qfac{-i + n} \qfac{-k + n} \qfac{1 + k + n}}
{\qfac{n} \qfac{-i - k + n} \qfac{1 - i + k + n}},
$$
we obtain
\begin{eqnarray*}
\sum_{j \ge 0}t^{(1)}_p(n, k, i, j) - \sum_{j \ge 0}t^{(2)}_p(n, k, i, j)  
&=& \sum_{j = 0}^{p-i}t^{(1)}_p(n, k, i, j) - 
\sum_{j = 0}^{p-i-1}t^{(2)}_p(n, k, i, j)  
\\
&=& a'_p(n, i) \frac{s_p(n-i, k)}{s_p(n, k)},
\end{eqnarray*}
and therefore
\begin{eqnarray*}
\sum_{i \ge 1}\sum_{j \ge 0}
\left(t^{(1)}_p(n, k, i, j)-t^{(2)}_p(n, k, i, j)\right)
&=& \sum_{i=1}^{p} a'_p(n, i) \frac{s_p(n-i, k)}{s_p(n, k)} \\
&=& r_p(n,k).
\end{eqnarray*}

We are going to show that 
\begin{equation}
\label{eqn.cert}
1 + r_p(n, k) 
= \Cert_p(n, k) - \Cert_p(n, k-1) \frac{s_p(n, k-1)}{s_p(n, k)}. 
\end{equation}

If \eqref{eqn.cert} is true, we can multiply both sides by $s_p(n, k)$ and
obtain
$$
s_p(n, k) + \sum_{i=1}^{p}a'_p(n, i) s_p(n-i, k)= 
\Cert_p(n, k) s_p(n, k) - \Cert_p(n, k-1) s_p(n, k-1).
$$
Summing over $k \ge 0$, and using telescoping summation of the right
hand side, and the boundary condition $\Cert_p(n, -1) = 0$, completes
the proof of \eqref{eqn.recursion.backwardshift}.
Notice incidentally that
$\Cert_p(n, k)$ is the corresponding certificate of
\eqref{eqn.recursion} in the WZ algorithm.

A recursion for both of the functions $t^{(h)}_p(n, k, i, j),\ h=1,2,$ is
\begin{equation}
\begin{array}{l}
- q^p (q^k - q^n) (-q + q^n) (-1 + q^{1 + k + n}) 
t^{(h)}_{p-1}(-1 + n, k, -1 + i, j) \\ [2.5mm]
+ q^{2 + k + 2 n} (-1 + q^n) t^{(h)}_{p-2}(n, k, i, -1 + j) \\[2.5mm]
- q^{2 + k + 2 n} (-1 + q^n) t^{(h)}_{p-1}(n, k, i, -1 + j) \\[2.5mm]
- q^{2 + k} (-1 + q^n) t^{(h)}_{p-1}(n, k, i, j)
+ q^{2 + k} (-1 + q^n) t^{(h)}_p(n, k, i, j) = 0.     
\end{array}     \nonumber
\end{equation}
This can be checked by dividing the equation by $t^{(h)}_p(n, k, i, j)$
and then both sides are rational functions in $q,q^n,q^p,q^k$; the identity
can then be checked easily.

Summing over $i \ge 1$ and $j \ge 0$, and noticing that 
$t^{(h)}_p(n, k, i, -1) = 0$,
we get
\begin{equation}
\begin{array}{l}
- q^p (q^k - q^n) (-q + q^n) (-1 + q^{1 + k + n}) r_{p-1}(n-1, k) \\[2.5mm]
+ q^{2 + k + 2 n} (-1 + q^n) r_{p-2}(n, k)     
- q^{2 + k + 2 n} (-1 + q^n) r_{p-1}(n, k))                      \\[2.5mm]
- q^{2 + k} (-1 + q^n) r_{p-1}(n, k)
+ q^{2 + k} (-1 + q^n) r_p(n, k)                                \\[2.5mm]
=q^p (q^k - q^n) (-q + q^n) (-1 + q^{1 + k + n})
\left(t^{(1)}_{p-1}(n, k, 0, 0) 
- t^{(2)}_{p-1}(n, k, 0, 0)\right)           \\[2.5mm]                     
= q^p (-q + q^n) (-q^k + q^n) (-1 + q^{1 + k + n}).\label{eqn.recur.ds2}
\end{array}
\end{equation}
What is left to prove now is that 
$D_p(n, k)$ satisfies the same recursion
as in \eqref{eqn.recur.ds2},
and \eqref{eqn.cert} is true 
for all $n$ when $p = 1$ and 2. 
Checking the former assertion is  
simple arithmetic since $D_p(n,k)$ is a rational function,
while the latter can be proved by checking \eqref{eqn.recursion.backwardshift} 
directly for $p = 1$ and 2.
For any specific $p$, $s_p(n, k)$ is hypergeometric, 
so this can be done using any of the software packages developed for the WZ 
algorithm; see for example \cite{PR3}.

When $p<0$, we can define\eject
\begin{eqnarray*}
t^{(1)}_p(n, k, i, j) &       =       &       
(-1)^{-i} q^{-i(2n-i+3)/2+(2p + i + 1)n+(2  n - i) j}   \\[2.5mm]
                        &               &
\frac{\qfac{n-1} \qfac{-k + n} \qfac{1 + k + n}}
{\qfac{n} \qfac{-i - k + n} \qfac{1 - i + k + n}}
\qbinom{-p - j - 1}{-p - i - j}\qbinom{i + j}{j},
\\[2.5mm]
t^{(2)}_p(n, k, i, j) &       =       &       
(-1)^{-i} q^{-i(2n-i+3)/2+(2p + i + 1)n+(2  n - i) j}   \\[2.5mm]
                        &               &
\frac{\qfac{n-1} \qfac{-k + n} \qfac{1 + k + n}}
{\qfac{n} \qfac{-i - k + n} \qfac{1 - i + k + n}}
\qbinom{-p - j - 2}{-p - i - j - 1}\qbinom{i + j}{j}, 
\end{eqnarray*}
and follow the same steps as above,  
where we only need to mention that both of the functions satisfy the 
following recursion
\begin{equation}
\begin{array}{l}
-q^p (q^k - q^n) (-q + q^n) (-1 + q^{1 + k + n}) t_{p-1}(n-1, k, i-1, j)   
 \\[2.5mm]
+ q^{2 + k + 2 n} (-1 + q^n) t_{p-2}(n, k, i, j)
- q^{2 + k + 2 n} (-1 + q^n) t_{p-1}(n, k, i, j-1)    \\[2.5mm]
- q^{2 + k} (-1 + q^n) t_{p-1}(n, k, i, j) 
+ q^{2 + k} (-1 + q^n) t_{p}(n, k, i, j-1) = 0.
\end{array}     \nonumber
\end{equation}
So far, we have shown that $\hatJ_p(n)$ is annihilated by an explicit
operator $C_p(E,Q,q)$:
$$
C_p(E,Q,q) \hatJ_p =0.
$$
If we prove that the above recursion has minimal $E$-degree, it will follow
that $C_p(E,Q,q)$ is indeed the non-commutative $C$--polynomial of the twist
knot $K_p$.
Since $C_p(E,Q,q)$ is monic in $E$, minimality will follows from the fact
that the polynomial $C_p(E,Q,1)$ is irreducible over $\BQ[E,Q]$.
This in turn follows from part (b) of \fullref{thm.CJp} and by
the fact that the $A$--polynomial of twist knots is irreducible 
(see \cite{HS2}).
This concludes the proof of \fullref{thm.1}.
\end{proof}

\section{Odds and ends}
\label{sec.que}

\subsection{Motivation for the rational map $\phi$}
\label{sub.motivation}

In this section we give some motivation for the strange-looking rational
map $\phi$. We warn the reader that this section is heuristic, and not
rigorous. However, it provides a good motivation.

Let us fix a sequence:
$$
f\co  \BN \longto \BQ(q)
$$
and let 
$$
g=\hat f\co  \BN \longto \BQ(q)
$$ 
be defined by:
$$
g(n)=\sum_{k=0}^\infty C(n,k) g(k),
$$
where $C(n,k)$ are as in Equation \eqref{eq.Cnk}. Let us suppose that
$f(k)$ is annihilated by an operator 
$$
P_f(E_k,Q_k,q)=\sum_{j=0}^d a_j(q,q^j) E_k^j.
$$ 

The question is to find
(at least heuristically) an operator $P_g(E,Q,q)$ that annihilates $g(n)$.
To achieve this, we will work in the Weyl algebra $\calW$
generated by the operators
$E,Q,E_k$ and $Q_k$ with the usual commutation relations.
\vspace{2pt}

Since $C(n,k)$ is closed form, a calculation shows that:
\begin{eqnarray*}
\frac{C(n+1,k)}{C(n,k)} &=& 
\frac{(1-q^{-n})(1-q^{1+k+n})}{(1-q^{k-n})(1-q^{1+n})} \\
\frac{C(n,k+1)}{C(n,k)} &=& -q^{-1-k}(1-q^{1+k-n})(1-q^{1+k+n}).
\end{eqnarray*}
In other words, $C(n,k)$ is annihilated by the left ideal in $\calW$ 
generated by the operators $P_1$ and $P_2$ where
\begin{eqnarray*}
P_2 &=& (1-Q^{-1})(1-qQQ_k) -(1-Q_k Q^{-1})(1-qQ) E \\
P_1 &=& -q^{-1}Q_k^{-1}(1-qQ_k Q^{-1})(1-qQ_k Q) -E_k
\end{eqnarray*}

\begin{lemma}
\label{lem.heuristic1}
$g(n)$ is annihilated by the operators $P_1$ and $P$ where
$$
P= \sum_{j=0}^d a_j(q,q^k) \frac{C(n,k+d)}{C(n,k+j)} E^j_k.
$$
\end{lemma}

\begin{proof}
It is easy to see that $g(n)$ is annihilated by $P_1$. Moreover,
\begin{eqnarray*}
0 & = & \sum_{j=0}^d a_j(q,q^k) f(k+j) \\
& = & \frac{1}{C(n,k+d)}  \sum_{j=0}^d a_j(q,q^k) 
\frac{C(n,k+d)}{ C(n,k+j)} C(n,k+j) f(k+j) \\
& = & \frac{1}{C(n,k+d)} P g(n).
\end{eqnarray*}
Thus, $P$ annihilates $g(n)$.
\end{proof}

According to Sister Celine's algorithm, we want to eliminate $Q_k$ (thus
obtaining $k$--free operators), and then set $E_k=1$. This will produce 
an operator in $E,Q$ and $q$ that annihilates $g(n)$. Finally, after setting
$q=1$, we will get a polynomial which contains the characteristic polynomial
of $g(n)$.

Now, here comes the heuristic: let us commute the evaluation at $q=1$
from last to first, and denote it by $\e$. 
Let us define two rational functions $R_1, R_2 \in \BQ(Q,Q_k)$ by:
\begin{eqnarray*}
R_1(Q,Q_k) &=& 
\frac{(1-Q^{-1})(1-Q_k Q)}{(1-Q_k Q^{-1})(1-Q)} \\
R_2(Q,Q_k) &=& -Q_k^{-1}(1-Q_k Q^{-1})(1-Q_k Q). 
\end{eqnarray*}
Observe that
$$
\e \frac{C(n,k+d)}{C(n,k+j)} = R_2(Q,Q_k)^{d-j}.
$$
Thus, by the above calculation,
\begin{eqnarray*}
\e P(E,Q,E_k,Q_k) &=& R_2^d \e P(R_2^{-1} E_k,Q_k) \\
\e P_1 &=& R_1(Q,Q_k)-E.
\end{eqnarray*}
Now, we want to eliminate $Q_k$ and then set $E_k=1$. The relation
$\e P_1=0$ is {\em linear} in $Q_k$. Solving, we obtain that:
$$
Q_k=\frac{1+QE}{Q+E}.
$$
Substituting this into $R_2^d \e P_f(R_2^{-1} E_k,Q_k)$ and setting $E_k=1$, 
we obtain:
\begin{align*}
R_2^d \e P_f(R_2^{-1} E_k,Q_k)&=R_3^d \e P(R_3, (1+QE)/(Q+E),1)\\&=
\e P^{\op}_f(R_3,(1+QE)/(Q+E),1),
\end{align*}
where
$$
R_3(E,Q)=\frac{E(Q^2-1)^2}{Q(E+Q)(1+EQ)}.
$$
In other words, after we rename $(E,Q)$ to $(L,M)$, we have:
\begin{eqnarray*}
\e P(L,M) &=& \e P^{\op}_f(R_3,(1+QE)/(Q+E),1) \\
&=& \phi \e P^{\op}(E,Q,1).
\end{eqnarray*}
In other words, we expect the characteristic polynomials $P_g(L,M)$
and $P_f(E,Q)$ of $g$ and $f$ to be related by:
\begin{equation}
\label{eq.Pfg}
P_g(L,M) \equiv_{L,M} \phi \, P^{\op}_f(E,Q),
\end{equation}
where $\equiv_{L,M}$ means equality, up to multiplication by monomials in $L$,
$M$, $L+M$, $1+LM$, $M-1$ and $M+1$. This is exactly how we came
up with the strange looking rational map $\phi$, and with \fullref{thm.CJp}.

In general, Equation \eqref{eq.Pfg} does not take into account repeated
factors in $P_g(L,M)$ and $P_f^{\op}(E,Q)$. Let us make this more precise.
Given $G \in \BQ[L,M]$, let us factor 
$G=u \prod_i G_i^{n_i}$ where $u$ is a unit, and $G_i$ are irreducible 
and $n_i \in \BZ$. This is possible, since $\BQ[L,M]$ is a unique 
factorization domain. Now, let us define 
$$
\rad(G)=\prod_i G_i^{\sgn(n_i)}
$$
to be the square-free part of $G$. 
\vspace{2pt}

Then, the equation
\begin{equation}
\label{eq.Pfg2}
\rad(G(L,M)) \equiv \rad(\phi \, F(E,Q)),
\end{equation}
implies that $F$ determines $G$ up to multiplication by suitable monomials,
and up to repeated factors. We may also invert the above equation,
keeping in mind that the map $\phi$ is 2-to-1.
\vspace{2pt}

Let us end this heuristic section with a lemma that sheds some light into
a possible relation between the $A$ and the $C$--polynomials of a knot.
\vspace{2pt}

Let $\equiv'$ denote equality of rational functions in $E,Q$ modulo
multiplication by monomials in $E, Q, 1+2Q+Q^2+QS$ and $1-2Q+Q^2+QS$.
These are precisely the monomials that map under $\phi$ to monomials in
$L,M,L+M,1+LM,M-1$ and $M+1$. 
\vspace{5pt}

\begin{lemma}
\label{lem.invertphi}
If $F$ and $G$ satisfy \eqref{eq.Pfg2},
then
\begin{align*}
\rad&(F(E,Q))\\& \equiv' 
\rad \left(
\mathrm{Res}_M\left(G\left(\frac{MQ-1}{M-Q} ,M\right),
E MQ-M^2 Q-Q+Q^2M+M \right) \right)
\end{align*}
where $\mathrm{Res}_M$ denotes the {\em resultant} with respect to $M$.
\end{lemma}
\vspace{5pt}

\begin{proof}
There is a geometric proof, which translates Equation \eqref{eq.Pfg2}
into the statement that $\phi$ induces a Zariski dense 
rational map of degree $2$:
$$
\{ (L,M) \in \BC^2 \, | \, G(L,M)=0 \} \longto
\{ (E,Q) \in \BC^2 \, | \, F(E,Q)=0 \}.
$$
Geometrically, it is clear that the domain determines the range and 
vice-versa.

There is an alternative algebraic proof. 
Let us try to invert the rational map $\phi$. In other words,
consider the system of equations
\begin{eqnarray*}
E &=& \frac{L(M^2-1)^2}{M(L+M)(1+LM)}; \\
Q &=& \frac{1+LM}{L+M},
\end{eqnarray*}
with $E,Q$ known and $L,M$ unknown. Solving the last with respect to $L$
gives:
$$
L=\frac{MQ-1}{M-Q}.
$$
Substituting into the first equation gives:
$$
E= M+M^{-1}-Q-Q^{-1}.
$$
Generically, this has two solutions in $Q$. Nevertheless, the above equation
is equivalent to
$$
E MQ-M^2 Q-Q+Q^2M+M=0.
$$
So, we can take resultant to eliminate $M$:
$$
\mathrm{Res}_M\left(G\left(\frac{MQ-1}{M-Q} ,M\right),
E MQ-M^2 Q-Q+Q^2M+M \right).
$$
The result follows.
\end{proof}

\subsection{Questions}
\label{sub.questions}

In this section we formulate several questions regarding the structure
and significance of the (non-commutative) $C$--polynomial of a knot.

Our first question may be thought of as a refined integrality property
for the cyclotomic function of a knot.

\begin{question}
\label{conj.1}
For which knots $K$, is the non-commutative $C$--polynomial monic
in $E$ with coefficients in $\BZ[Q^{\pm},q^{\pm}]$?
\end{question}

Motivated by \fullref{cor.Cphi}, it is tempting to formulate the 
following 

\begin{question}
\label{conj.2}
For which knots $K$, does the map $\phi$ of Equation \eqref{eq.phi}
give a Zariski dense rational map of degree $2$:
\begin{equation}
\label{eq.ACK}
\phi\co  \{(L,M) \in \BC^2 \, |\, A_K(L,M^{1/2})=0 \}
\longto \{(E,Q) \in \BC^2 \, |\, C_K^{\op}(E,Q)=0 \}.
\end{equation}
\end{question}

In view of \fullref{lem.invertphi}, \eqref{eq.ACK} is equivalent
to
\begin{equation}
\label{eq.conj2}
\rad(\phi C_K^{\op}(E,Q)) \equiv \rad(A(L,M^{1/2})).
\end{equation}

\begin{question}
\label{conj.3}
Is the genus of the $C$--polynomial $C_K(E,Q)$ of a knot always zero?
\end{question}

\begin{remark}
\label{rem.conj1}
It seems that Questions \ref{conj.1} and \ref{conj.3} cannot be positive
the same time for the $2$--bridge knot $K_{13/5}$
and for the simplest hyperbolic non-2--bridge knot $m082$.
\end{remark}

\begin{theorem}
\label{thm.CD}
If $C_K(E,1) \neq 0$, then the Alexander polynomial $\D_K(M)$ divides 
$$
C^{\op}_K(M-2+M^{-1},1).
$$
\end{theorem}

The above theorem provides a nontrivial consistency check of \fullref{conj.2}. Indeed, in \cite[Section 6]{CCGLS} Cooper et al. 
prove that 
$A_K(1,M^{1/2})$ is divisible by the Alexander polynomial $\D_K(M)$
at least when the latter has unequal complex roots.
On the other hand
$$
\phi(E)|_{L=1}=M-2+M^{-1}, \qquad \phi(Q)|_{L=1}=1.
$$
Thus, 
$$
\phi C^{\op}_K(E,Q)|_{L=1}=C^{\op}_K(M-2+M^{-1},1).
$$
Thus, if a knot satisfies \fullref{conj.2}, then $\D_K(M)$
divides $C^{\op}_K(M-2+M^{-1},1)$. This is precisely \fullref{thm.CD}.

\begin{question}
\label{que.1}
Does the $C$--polynomial of a knot 
have a classical geometric definition? 
\end{question}

In other words, we are asking for a {\em geometric meaning} of the 
rational map $\phi$ of Equation \eqref{eq.phi}.

\begin{question}
\label{que.2}
Is there any relation between the bi-graded knot invariant $C_K(E,Q,1)$
and some version of Knot Floer Homology? \fullref{thm.CD} states that
under mild hypothesis, $C_K(M-2+M^{-1},1,1)$ is divisible by the 
Alexander polynomial of $K$.
\end{question}

\subsection[Proof of \ref{thm.CD}]{Proof of \fullref{thm.CD}}
\label{sec.thm12}

\begin{proof}
The proof utilizes the algebra of generating functions and 
the fact that the generating function of the cyclotomic
function of a knot (evaluated at $1$) is given by the inverse Alexander 
polynomial. For a reference of the latter statement, see \cite{GL2}.

Now, let us give the details of the proof. 
We start from the recursion relation of the cyclotomic function:
$$
C_K(E,Q,q) \hatJ_K=0.
$$
Let us evaluate at $q=1$, and set
$$
C_K(E,1,1)=\sum_{j=0}^d a_j E^j, \qquad
I_K(n)=\hatJ_K(n)|_{q=1}.
$$
Then, we have for all $n$:
\begin{equation}
\label{eq.genf}
\sum_{j=0}^d a_j I_K(n+j)=0.
\end{equation}
Let us use the generating function:
$$
F_K(z)=\sum_{n=0}^\infty I_K(n) z^n.
$$
Equation \eqref{eq.genf} implies that
\begin{eqnarray*}
0 & = & 
\sum_{n=0}^\infty \sum_{j=0}^d a_j I_K(n+j) z^n \\
& = & \sum_{j=0}^d a_j z^{-j} \sum_{n=0}^\infty I_K(n+j) z^{n+j} \\
& = & \sum_{j=0}^d a_j z^{-j} F_K(z) + \quad \text{terms}(z) \\
& = & z^{-d} C^{\op}_K(z,1,1) F_K(z) + \quad \text{terms}(z),
\end{eqnarray*}
where $\text{terms}(z)$ is a Laurent polynomial in $z$.
Thus, assuming that $C^{\op}(z,1,1) \neq 0$, it follows that
$$
F_K(z) = -\frac{\text{terms}(z)}{C^{\op}_K(z,1,1)}.
$$
The Melvin--Morton--Rozansky Conjecture (proven by Bar-Natan and the first 
author in \cite{B-NG}), together with the cyclotomic expansion of the colored
Jones function implies that 
$$
F_K(M-2+M^{-1})=\frac{1}{\D_K(M)}.
$$
For example, see \cite[Lemma 2.1]{GL2}. Thus, 
$$
-\frac{\text{terms}(M-2+M^{-1})}{C^{\op}_K(M-2+M^{-1},1,1)}
=\frac{1}{\D_K(M)}.
$$
The result follows.
\end{proof}

\appendix

\section{A table of non-commutative $C$--polynomials}
\label{sec.table}

We finish with a table of the non-commutative $C$--polynomial
of twist knots $K_p$ for $p=-3,\dots,3$, taken from \fullref{thm.1}.
In each matrix, the upper left entry indicates the $C_p(E,Q,q)$ polynomial
and the entries in the $E^i$--row and $Q^j$--column indicate the coefficient
of $Q^j E^i$ in  $C_p(E,Q,q)$. For example, $C_1(E,Q)=E+q^2Q$.
$$
\mat
{
C_1 & Q^0 & Q^1 \\      \hline
E^0 & 0 & q^2 \\
E^1 & 1 & 0 
}
\hspace{2cm}
\mat
{
C_2 & Q^0 & Q^1 & Q^2 & Q^3  \\ \hline
E^0 & 0 & 0 & q^6 & -q^7 \\
E^1 & 0 & q^3+q^4 & -q^5 & q^7 \\
E^2 & 1 & 0 & 0 & 0 
}
$$
{\small$$
\mat
{
C_2 & Q^0 & Q^1 & Q^2 & Q^3 & Q^4 & Q^5  \\     \hline
E^0 & 0 & 0 & 0 & q^{12} & -q^{13}-q^{14} & q^{15} \\
E^1 & 0 & 0 & q^8 + q^9 + q^{10} & -q^{10}-2 q^{11} -q^{12} & 2 q^{13}+q^{14}
& -q^{15}-q^{16} \\
E^2 & 0 & q^4 +q^5+q^6 & -q^7-q^8 & q^{10}+q^{11} & -q^{13} & q^{16} \\
E^3 & 1 & 0 & 0 & 0 & 0 & 0
}
$$}
$$
\mat{
C_{-1} & Q^1 \\ \hline 
E^0  & -1 \\
E^1 & 1
}
\hspace{2cm}
\mat{
C_{-2} & Q^{-2} & Q^{-1} & Q^0 \\       \hline 
E^0 & 0 & -q^{-2} & q^{-1} \\
E^1 & -q^{-4} & q^{-2} & -q^{-1}-1 \\
E^2 & 1 & 0 & 0
}
$$
$$
\mat{
C_{-3}  & Q^{-4} & Q^{-3} & Q^{-2} & Q^{-1} & Q^0 \\    \hline 
E^0 &  0 & 0 & -q^{-6} & q^{-5} + q^{-4} & -q^{-3} \\
E^1 &  0 & -q^{-9} -q^{-8} & q^{-7} + 2 q^{-6} & -q^{-5}-2 q^{-4}-q^{-3} &
q^{-3}+q^{-2}+q^{-1} \\
E^2 & -q^{-12} & q^{-9} & -q^{-7} -q^{-6} & q^{-4}+q^{-3} & -q^{-2}-q^{-1}-1\\
E^3 & 1 & 0 & 0 & 0 & 0
}
$$

\bibliographystyle{gtart}
\bibliography{link}

\end{document}